\documentclass{amsart} 
\usepackage{mathrsfs,amscd,amssymb,graphicx,latexsym,xypic}
\newcommand{\LATTICES}{\mathscr{L}}

\newcommand{\POSPACES}{\mathscr{P}}
\newcommand{\PREORDEREDSPACES}{\mathscr{Q}}

\newcommand{\SPACES}{\mathscr{T}}
\newcommand{\id}{\mathrm{id}}

\newcommand{\upper}[2]{\leqslant_{#1}\!\![#2]}
\newcommand{\downer}[2]{\leqslant_{#1}^{-1}\!\![#2]}

\newtheorem{theorem}{Theorem}[section]
\newtheorem{lemma}[theorem]{Lemma}
\newtheorem{proposition}[theorem]{Proposition}

\theoremstyle{definition}
\newtheorem{definition}[theorem]{Definition}
\newtheorem{example}[theorem]{Example}

\begin{document}                                                                
\title{Criteria for homotopic maps to be so along monotone homotopies}
\author{Sanjeevi Krishnan}
\address{Laboratoire d'Informatique de l'\'{E}cole Polytechnique\\Palaiseau, France}
\begin{abstract}
	The state spaces of machines admit the structure of time.
	A homotopy theory respecting this additional structure can detect machine behavior unseen by classical homotopy theory.
	In an attempt to bootstrap classical tools into the world of abstract spacetime, we identify criteria for classically homotopic, monotone maps of pospaces to \textit{future homotope}, or homotope along homotopies monotone in both coordinates, to a common map.
	We show that consequently, a hypercontinuous lattice equipped with its Lawson topology is \textit{future contractible}, or contractible along a future homotopy, if its underlying space has connected CW type.
\end{abstract}
\maketitle
\section{Introduction}

\begin{figure}[h!]
	\begin{center}
		\includegraphics[width=40mm,height=40mm]{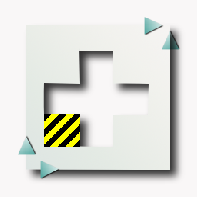} 
	\end{center}
	\caption{State space of a binary semaphore, as in \cite[Figure 7]{fgr:ditop}}
	\label{fig:2sem}
\end{figure}

The state spaces of machines often admit partial orders which describe the causal relationship between states.  
For example, the unit interval $\mathbb{I}$ equipped with its standard total order represents the states of a finite, sequential process.
Figure \ref{fig:2sem} illustrates the state space $X$ of two sequential processes accessing a binary semaphore.
Thinking of the upper corner as the desired end state, we view monotone paths $\mathbb{I}\rightarrow X$ reaching the striped zone as unsafe executions of our binary system, doomed never to terminate successfully.
We can thus articulate critical machine behavior in the language of partially ordered spaces.

\begin{figure}[h!]
	\begin{center}
		\begin{tabular}{cc}
			\includegraphics[width=40mm,height=40mm]{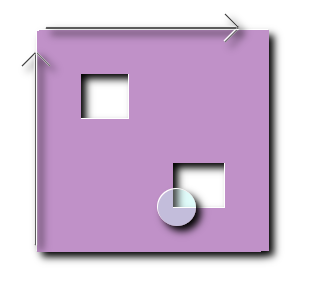} 
			& \includegraphics[width=40mm,height=40mm]{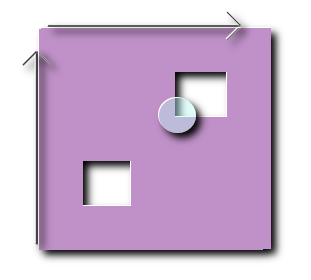} \\
			(a)
			& (b) 
		\end{tabular}
	\end{center}
	\caption{Partially ordered state spaces, as in \cite[Figure 14]{fgr:ditop}}
	\label{fig:inequivalent}
\end{figure}

A homotopy theory respecting this additional structure of time potentially can detect machine behavior invisible to classical homotopy theory, as demonstrated in \cite{fgr:ditop}.
A suitable theory should distinguish between the homotopy equivalent state spaces given in Figure \ref{fig:inequivalent}, for example.
In an attempt to exploit classical arguments in a homotopy theory of preordered spaces, we seek criteria under which two homotopic, monotone maps $X\rightarrow Y$ of pospaces are in fact homotopic through monotone maps.
Certain cubical approximation results in \cite{fajstrup:approx} implicitly use one such criterion: when $Y$ is a convex sub-pospace of an ordered topological vector space.
Lemma \ref{lem:main.result} identifies alternative criteria which do not require vector space structures: when $X$ is a compact pospace whose ``lower'' sets generated by open subsets are open and $Y$ is a continuous lattice equipped with its Lawson topology.

We can further refine classical homotopy theory, following \cite{grandis:d}.
Consider two monotone maps $f,g:X\rightarrow Y$ of preordered spaces.  
We say that $f$ \textit{future homotopes to} $g$ if a classical homotopy from $f$ to $g$ defines a monotone map $X\times\mathbb{I}\rightarrow Y$.
We call a preordered space \textit{future contractible} if the identity on it future homotopes to a constant map.  
Lemma \ref{lem:future.homotope} identifies criteria under which two monotone maps $X\rightarrow Y$ homotopic through monotone maps future homotope to a common map: when $X$ is compact Hausdorff and $Y$ is the order-theoretic dual of a continuous lattice $L$ equipped with the dual Lawson topology of $L$.
We obtain the following consequence.

\vspace{2mm}
\noindent {\bf Proposition \ref{prop:dicontract}}
\hspace{1mm}\textit{A hypercontinuous lattice equipped with its Lawson topology is future contractible if its underlying space has connected CW type.}\\
\vspace{2mm}

It follows that a hypercontinuous lattice equipped with its Lawson topology is ``past'' contractible if its underlying space has connected CW type, by symmetry.
In \S\ref{sec:pospaces}, we review some basic definitions, examples, and properties of preordered spaces.
In \S\ref{sec:homotopies} we prove Lemmas \ref{lem:main.result} and \ref{lem:future.homotope}, followed by Proposition \ref{prop:dicontract}.

\section{Preordered spaces}\label{sec:pospaces}

A \textit{preordered space} is a preordered set equipped with a topology.
An example is a \textit{topological sup-semilattice} (\textit{inf-semilattice}), a sup-semilattice (inf-semilattice) equipped with a topology making the binary $\sup$ ($\inf$) operator jointly continuous.  
A \textit{monotone map} is a continuous, (weakly) monotone function between preordered spaces.
The forgetful functor
$$U:\PREORDEREDSPACES\rightarrow\SPACES$$
from the category $\PREORDEREDSPACES$ of preordered sets and monotone functions to the category $\SPACES$ of spaces and continuous functions has a left adjoint.
We write $\ddot{U}:\PREORDEREDSPACES\rightarrow\PREORDEREDSPACES$ for the composite of $U$ with its left adjoint, and we write $\epsilon:\ddot{U}\rightarrow\id_{\PREORDEREDSPACES}$ for the counit of the adjunction.

For each preordered space $X$, we write $\leqslant_X$ for its preorder and
$$\upper{X}{A}=\bigcup_{a\in A}\{x\;|\;a\leqslant_Xx\},\quad\downer{X}{A}=\bigcup_{a\in A}\{x\;|\;x\leqslant_Xa\}$$
for the ``upper'' and ``lower'' sets, respectively, generated by a subset $A\subset X$.

\begin{figure}
	\begin{center}
		\begin{tabular}{cc}
			\includegraphics[width=30mm,height=30mm]{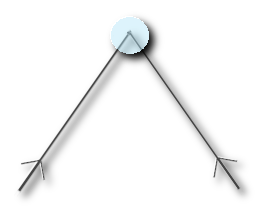} 
			& \includegraphics[width=30mm,height=30mm]{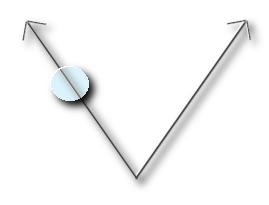} \\
			$X_1$ & $X_2$
		\end{tabular}
	\end{center}
	\caption{Compact pospaces with and without lower open partial orders.}
	\label{fig:lower.open}
\end{figure}

\begin{example}
	In Figure \ref{fig:lower.open}, $X_1$ is a topological sup-semilattice and
	$$\downer{X_1}{V_1}=X_1$$
	for $V_1$ the circled open subset of $X_1$.
\end{example}

\begin{example}
	In Figure \ref{fig:lower.open}, $X_2$ is a topological inf-semilattice and
	$$\downer{X_2}{V_2}$$
	is not open in $X_2$, for $V_2$ the circled open subset of $X_2$.
\end{example}

\begin{example}[Counterexamples]
	The pospaces of Figure \ref{fig:lower.open} are neither inf-semilattices nor sup-semilattices, even though their underlying posets are complete lattices.
\end{example}

Certain preorders are ``continuous'' in the following sense.

\begin{definition}
	A preorder $\leqslant_X$ on (the points of a) space $X$ is \textit{lower open} if 
	$$\downer{X}{V}$$
	is open in $X$ for each open subset $V\subset X$.
\end{definition}

An example of a lower open preorder is the trivial preorder on a space.  
The class of preordered spaces having lower open preorders is closed under products and coproducts.

\begin{lemma}\label{lem:non.branching}
	All topological sup-semilattices have lower open preorders.
\end{lemma}
\begin{proof}
	For each open subset $V$ of a topological $\sup$-semilattice $L$, 
	$$\downer{L}{V}=\pi_2((V\times L)\cap\sup\!^{-1}(V)),$$
	where $\pi_2:L\times L\rightarrow L$ denotes projection onto the second factor, is open in $L$ because $\pi_2$ is an open map and $\sup$ is a continuous function $L\times L\rightarrow L$.
\end{proof}

Recall from \cite{scott:lattices} that a \textit{pospace} is a preordered space $X$ whose partial order $\leqslant_X$ is antisymmetric ($x\leqslant_Xy\leqslant_Xx$ implies $x=y$) and has closed graph in the standard product topology $X\times X$.  

\begin{example}
	The preordered spaces in all of the figures are pospaces.
\end{example}

Pospaces are automatically Hausdorff by \cite[Proposition VI-1.4]{scott:lattices}.
Examples include Hausdorff topological sup-semilattices and Hausdorff topological inf-semilattices by \cite[Proposition VI-1.14]{scott:lattices}.
In particular, \textit{continuous lattices equipped with their Lawson topologies}, which \cite[Theorem VI-3.4]{scott:lattices} characterizes as compact Hausdorff, topological inf-semilattices which have a maximum and whose points admit neighborhood bases of sub-semilattices, are pospaces.

We can construct the ``free continuous lattice generated by a compact pospace,'' following \cite[Example VI-3.10 (ii)]{scott:lattices}.
Let $\POSPACES$ denote the full subcategory of $\PREORDEREDSPACES$ consisting of compact pospaces.
Inclusion $i:\LATTICES\hookrightarrow\POSPACES$ from the category $\LATTICES$ of continuous lattices equipped with their Lawson topologies and continuous semilattice homomorphisms preserving maxima has a left adjoint
$$F:\POSPACES\rightarrow\LATTICES$$
sending each compact pospace $X$ with topology $\mathcal{T}_X$ to the poset of all closed subsets $C\subset X$ satisfying $C=\upper{X}{C}$, ordered by reverse inclusion and having topology generated by the subsets
$$\{A\;|\;A\subset V\},\;\;\{B\;|\;B\cap W\neq\varnothing\},\quad V,W\in\mathcal{T}_X,\;W=\downer{X}{W}.$$

The unit is the natural map $\upsilon_X:X\rightarrow FX$ defined by $x\mapsto\upper{X}{\{x\}}$.  
The counit is the infinitary infimum operator $\bigwedge:FL\rightarrow L$.

\begin{lemma}\label{lem:alternative.criterion}
	Consider a compact pospace $X$. 
	The inclusion
	$$FX\hookrightarrow FUX$$
	is continuous if $\leqslant_X$ is lower open.
\end{lemma}
\begin{proof}
	Consider an open subset $W\subset X$.
	The set
	\begin{eqnarray*}
		\{B\in FX\;|\;B\cap W\neq\varnothing\}
		&=& \{B\in FX\;|\;\upper{X}{B}\cap W\neq\varnothing\}\\
		&=& \{B\in FX\;|\;B\;\cap\downer{X}{W}\neq\varnothing\}
	\end{eqnarray*}
	is open in $FX$ if $\downer{X}{W}$ is open in $X$.  
	The claim then follows.  
\end{proof}

We can thus give a useful recipe for converting continuous functions into monotone maps. 

\begin{lemma}\label{lem:retraction}
	For each compact pospace $X$ having lower open partial order and each continuous lattice $Y$ equipped with its Lawson topology, the function
	$$U:\POSPACES(X,Y)\rightarrow\SPACES(UX,UY)$$
	has a retraction $f\mapsto(x\mapsto\bigwedge f(\leqslant_X[\{x\}]))$.
\end{lemma}
\begin{proof}
	For a continuous function $f:UX\rightarrow UY$, the composite function
	$$\xymatrix@C35pt{
	X\times\mathbb{I}\ar[r]^{\!\!\!\!\!\!\upsilon_{X\times\mathbb{I}}} 
	& F(X\times\mathbb{I})\ar[r]^{\!\!j} 
	& F\ddot{U}(X\times\mathbb{I})\ar[r]^{\;\;\;\;\;F(f)}
	& F\ddot{U}Y\ar[r]^{\;F(\epsilon_Y)}
	& FY\ar[r]^{\bigwedge}
	& Y,
	}$$
	where $j$ denotes the inclusion function, is a monotone map by Lemma \ref{lem:alternative.criterion}.
	This composite sends $x$ to $\bigwedge f(\leqslant_X[\{x\}])$, which equals $f(x)$ if $f$ is monotone.
\end{proof}

\section{The homotopy theory}\label{sec:homotopies}

We refine the classical homotopy relation, first by defining the ``dihomotopy'' relation of \cite{fgr:ditop}.
Let $\mathbb{I}$ be the unit interval $[0,1]$ equipped with its standard total order.  
Fix preordered spaces $X,Y$.  
For every pair of monotone maps
$$f,g:X\rightarrow Y,$$
we write $f\sim g$ if $f$ is homotopic through monotone maps to $g$, or equivalently, if a homotopy $Uf\sim Ug$ defines a monotone map $X\times\ddot{U}\mathbb{I}\rightarrow Y$.
Following classical notation, let $[f]$ denote the $\sim$-class of a monotone map $f:X\rightarrow Y$, and let $[X,Y]$ denote the set of all such equivalence classes $[f]$.
The forgetful functor $U:\POSPACES\rightarrow\SPACES$ to the category $\SPACES$ of spaces induces a natural function 
\begin{equation}
	U_*:[X,Y]\rightarrow[UX,UY]
	\label{eqn:group.completion}
\end{equation}
to the set of homotopy classes $[UX\rightarrow UY]$ of continuous functions $UX\rightarrow UY$.

\begin{figure}
	\begin{center}
		\begin{tabular}{cc}
			\includegraphics[width=30mm,height=30mm]{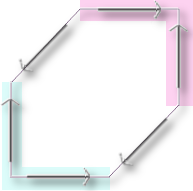} 
			& \includegraphics[width=30mm,height=30mm]{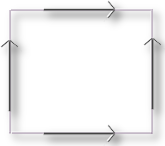} \\
			$X_3$ & $X_4$ 
		\end{tabular}
	\end{center}
	\caption{$U_*:[X_3,X_4]\rightarrow[UX_3,UX_4]$ neither injective nor surjective.}
	\label{fig:not.dihomotopic}
\end{figure}

\begin{example}\label{eg:non.injective}
	Consider the pospaces given in Figure \ref{fig:not.dihomotopic}.  
	The monotone map $X_3\rightarrow X_4$ surjectively wrapping the lower blue corner around $X_4$ is homotopic, though not through monotone maps, to a monotone map $X_3\rightarrow X_4$ surjectively wrapping the upper red corner around $X_4$. 
	Thus (\ref{eqn:group.completion}) need not be injective. 
	No monotone map $X_3\rightarrow X_4$ has Brouwer degree greater than $1$.
	Thus (\ref{eqn:group.completion}) need not be surjective. 
\end{example}

Directed homotopy theory reduces to classical homotopy theory and order-theory precisely when (\ref{eqn:group.completion}) is injective.  
The following lemma gives us such a case.

\begin{lemma}\label{lem:main.result}
	For each compact pospace $X$ having lower open partial order and each continuous lattice $Y$ equipped with its Lawson topology, the function
	$$U_*:[X,Y]\rightarrow[UX,UY]$$
	has a well-defined retraction $[f]\mapsto[x\mapsto\bigwedge f(\leqslant_X[\{x\}])].$
\end{lemma}
\begin{proof}
	For a compact pospace $A$ such that $\leqslant_A$ is lower open and a continuous lattice $B$ equipped with its Lawson topology, let $R_{A,B}:\SPACES(UA,UB)\rightarrow\POSPACES(A,B)$ denote the retraction defined by Lemma \ref{lem:retraction}.  
	The diagram
	$$
	\xymatrix@C50pt{
	\SPACES(UX\coprod UX,UY)\ar[r]^{\quad R_{(X\coprod X),Y}}\ar[d]_{(x\mapsto(x,0))\coprod(x\mapsto(x,1))}
	& \POSPACES(X\coprod X,Y)\ar[d]^{(x\mapsto(x,0))\coprod(x\mapsto(x,1))}\\
	\SPACES(UX\times U\mathbb{I},UY)\ar[r]_{\quad R_{X\times\ddot{U}\mathbb{I},Y}}
	& \POSPACES(X\times\ddot{U}\mathbb{I},Y),\\
	}
	$$
	is commutative and thus $R_{X,Y}$ passes to $\sim$-classes to define our desired retraction.
\end{proof}

\begin{figure}[h!]
	\begin{center}
		\begin{tabular}{cc}
			\includegraphics[width=40mm,height=40mm]{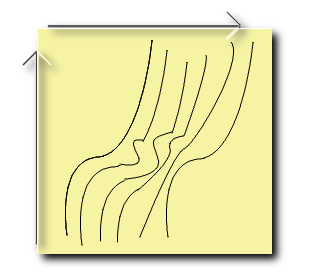} 
			& \includegraphics[width=40mm,height=40mm]{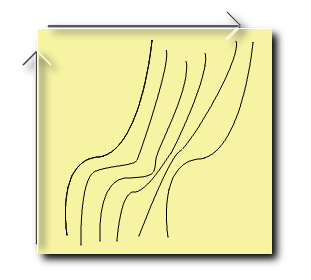} \\
			(a)
			& (b) 
		\end{tabular}
	\end{center}
	\caption{A classical homotopy (a) and a homotopy (b), obtained from an application of Lemma \ref{lem:retraction} to (a), through monotone paths}
	\label{fig:monotonization}
\end{figure}

\begin{example}
	Consider Figure \ref{fig:monotonization}.
	Under the retraction given in Lemma \ref{lem:retraction}, the homotopy through monotone paths in (b) is the image of the classical homotopy of paths in (a).
\end{example}

We refine the dihomotopy relation of \cite{fgr:ditop}, following \cite{grandis:d}.

\begin{definition}
	Given preordered spaces $X,Y$ and monotone maps
	$$f,g:X\rightarrow Y,$$
	we say that $f$ \textit{future homotopes to} $g$ if there exists a monotone map $h:X\times\mathbb{I}\rightarrow Y$ such that $h(-,0)=f$ and $h(-,1)=g$.
	A preordered space $X$ is \textit{future contractible} if $\id_X:X\rightarrow X$ future homotopes to a constant map.
\end{definition}

\begin{lemma}\label{lem:future.homotope}
	Consider a pair of monotone maps
	$$g_1,g_2:X\rightarrow Y$$
	from a compact Hausdorff preordered space $X$ to a Lawson semilattice $Y$, homotopic through monotone maps.  
	There exists a monotone map which future homotopes to both $g_1$ and $g_2$.
\end{lemma}
\begin{proof}
	Let $h:g_1\sim g_2$ be a homotopy through monotone maps.  
	The rules 
	$$j(x,t)=\bigwedge h(x,[0,1-t]),\quad k(x,t)=\bigwedge h(x,[t,1])$$
	define functions $j,k:X\times\mathbb{I}\rightarrow Y$.  
	The functions $j,k$ are continuous by Lemma \ref{lem:retraction} because $\leqslant_{\ddot{U}X\times\mathbb{I}}$ and its order-theoretic dual are lower open.  
	The functions $j,k$ are monotone because $\bigwedge$ is a monotone operator.  
	Thus $j(-,0)=k(-,0)$ future homotopes to $j(-,1)=h(-,0)=g_1$ and $k(-,1)=h(-,1)=g_2$.
\end{proof}

\begin{example}
	On hom-sets $\POSPACES(X,Y)$ for which $Y$ is a continuous lattice equipped with its Lawson topology, the dihomotopy relation $\sim$ coincides with the \textit{d-homotopy} relation of \cite{grandis:d}, as a consequence of Lemma \ref{lem:future.homotope}.
\end{example}

Recall that a space has \textit{connected CW type} if it is homotopy equivalent to a connected CW complex.
Recall from \cite{scott:lattices} that a hypercontinuous lattice is a continuous lattice whose Lawson and dual Lawson topologies agree.  
Thus a hypercontinuous lattice equipped with its Lawson topology is precisely a compact Hausdorff (inf- and sup-) topological lattice whose points admit, with respect to each semilattice operation, neighborhood bases of sub-semilattices.

\begin{proposition}\label{prop:dicontract}
	A hypercontinuous lattice equipped with its Lawson topology is future contractible if its underlying space has connected CW type.
\end{proposition}
\begin{proof}
	Consider a hypercontinuous lattice $L$ equipped with its Lawson topology, and suppose $UL$ has connected CW type.
	The space $UL$ is therefore path-connected.  
	Moreover, $UL$ has trivial homotopy groups because the binary $\inf$ operator gives $UL$ the structure of an associative, idempotent $H$-space.
	The map $\id_L$ is homotopic through monotone maps to a constant map $c$ taking the value $\max L$ by Lemma \ref{lem:main.result} - $\id_{UL}$ is homotopic to $U(c)$ by the Whitehead Theorem, $L$ is a compact pospace, and $\leqslant_L$ is lower open by Lemma \ref{lem:non.branching}.
	The map $\id_L$ and $c$ future homotope to $c$ by Lemma \ref{lem:future.homotope} because $L$ is the dual of a continuous lattice equipped with the dual Lawson topology of $L$.
\end{proof}

\section{Conclusion}
The state spaces of machines in nature arise as ``locally partially ordered'' geometric realizations of cubical complexes, as in \cite{fgr:ditop}.  
Such ``locally partially ordered'' spaces are hypercontinuous lattices precisely when they are continuous lattices, the computational steps of computable partially recursive functions in \cite{scott:outline}. 
Thus Proposition \ref{prop:dicontract} and Example \ref{eg:non.injective} suggest that the directed homotopy theories of \cite{fgr:ditop,grandis:d} measure at least some of the failure, undetected by classical homotopy theory, of a state space to represent a deterministic, computable process.

\section{Acknowledgements}
The author thanks Eric Goubault and Emmanuel Haucourt for their helpful comments and suggestions.
The author also thanks Zack Apoian, Steven Paulikas, and John Posch.

\end{document}